\input amstex
\documentstyle{amsppt}
\TagsOnRight
\overfullrule = 0pt

\vsize=9truein
\hsize=6.5truein
\hoffset=0pt
\voffset=-.5truein

\input pictex

\def\vp{\varepsilon}
\def\cl#1{{\Cal{#1}}}
\def\bb#1{{\Bbb{#1}}}



\def\HH{{\Bbb H}}
\def\KK{{\Bbb K}}

\def\QQ{{\Bbb Q}}
\def\RR{{\Bbb R}}

\def\ZZ{{\Bbb Z}}

\def\cK{{\Cal K}}

\def\D{{\Delta}}

\def\G{{\Gamma}}



\def\PGL{{\hbox{\rm{PGL}}}}
\def\SL{{\hbox{\rm{SL}}}}

\def\tto{{\hbox{\rm{II}}}_{1}}
\def\aut{{\hbox{\rm Aut}}}

\def\pgl{{\hbox{\rm{PGL}}}}
\def\psl{{\hbox{\rm{PSL}}}}

\def\SO{{\hbox{\rm{SO}}}}

\def\VN{{\hbox{\rm VN}}}

\def\hd{{\hbox{\rm hd}}}

\def\bs{\backslash}
\def\diam{{\hbox{\rm diam}}}

\topmatter
\title
STRONG SINGULARITY FOR SUBALGEBRAS \\ OF FINITE FACTORS
\endtitle

\author
Guyan Robertson\qquad Allan M.~Sinclair\qquad Roger R.~Smith
\endauthor
\address
Department of Mathematics,
University of Newcastle,
Callaghan, NSW 2308,
Australia\endaddress
\email
guyan\@frey.newcastle.edu.au
\endemail
\address
Department of Mathematics,
University of Edinburgh,
Edinburgh EH9 3JZ,
Scotland\endaddress
\email
allan\@maths.ed.ac.uk
\endemail
\address
Department of Mathematics,
Texas A\&M University,
College Station, TX 77843,
USA
\endaddress
\email
rsmith\@math.tamu.edu
\endemail
\thanks 2000 {\it Mathematics Subject Classification}. 46L10, 22D25.
\endthanks
\thanks
The first author was partially supported by the Australian Research
Council.
\endthanks
\thanks
The third author was partially supported by the National Science
Foundation.
\endthanks
\abstract
In this paper we develop the theory of strongly singular subalgebras
of von
Neumann algebras, begun in earlier work. We mainly examine the situation
of type $\tto$ factors arising from countable discrete groups. We give
simple criteria for strong singularity, and use them to construct strongly
singular subalgebras. We particularly focus on groups which act on
geometric objects, where the underlying geometry leads to strong
singularity.
\endabstract
\endtopmatter

\document

\head 1. Introduction\endhead

Let $\Cal A$ be a maximal abelian self--adjoint subalgebra (masa) in a
type $\tto$ factor $\Cal M$ with trace $tr$. In [6], Dixmier
identified various classes of masas based on the structure of the
normalizer
$$
N(\Cal A)=\{u\in \Cal M : u{\Cal{A}}u^*=\Cal A,
u{\hbox{ unitary}}\}.\tag{1.1}
$$
In particular, $\Cal A$ is said to be singular if $N(\Cal A)
\subseteq \Cal A$, so that the only normalizing unitaries already
belong to $\Cal A$. He also provided some examples of singular masas
inside factors arising from discrete groups. However, it is a difficult
problem to decide whether a given masa is singular, and this prompted the
second and third authors to introduce the concept of strong singularity in
[20]. The trace induces a norm $\|\cdot \|_2$ on $\Cal M$ by
$\|x\|_2=(tr(x^*x))^{1/2}$, and a norm $\|\cdot \|_{\infty,2}$ may then be
defined for a map $\phi : \Cal M \to \Cal M$ by
$$
\|\phi \|_{\infty,2}=\sup\{\|\phi(x)\|_2:\|x\|\leq 1\}.\tag{1.2}
$$
Letting ${\Bbb{E}}_{\Cal N}$ denote the unique trace preserving
conditional expectation onto any von Neumann subalgebra $\Cal N$ of
$\Cal M$, strong singularity of a masa $\Cal A$ (or of a general
von Neumann subalgebra) is then defined by requiring the inequality
$$
\|{\Bbb{E}}_{\Cal A}-{\Bbb{E}}_{u{\Cal A}u^*}\|_{\infty,2}
\geq \|(I-{\Bbb{E}}_{\Cal A})(u)\|_2\tag{1.3}
$$
to hold for all unitaries $u \in \Cal M$. Singularity of any such masa
is an immediate consequence of (1.3).
The objective in introducing this concept was to have an easily verifiable
criterion for singularity. For example, the masa in $VN({\Bbb{F}}_2)$
generated by one of the generators of ${\Bbb{F}}_2$ satisfies
(1.3), showing singularity (which was, of course, known to
Dixmier, [6]). However, in [20], the problem of exhibiting large
classes of strongly singular masas was not addressed beyond examples
arising from hyperbolic groups. The purpose of this paper is to examine
several general contexts in which strongly singular masas and subalgebras
appear naturally and in profusion. Since type $\tto$ factors are closely
connected to discrete groups, much of the work (but not exclusively) will
be in this area. We now give a brief outline of the paper.

The second section gives a criterion for determining when a von Neumann
subalgebra is strongly singular (Lemma 2.1). This is a minor
modification of a result from [20], but the new version is slightly
more flexible and thus more widely applicable. We use it to generate
classes of strongly singular masas based on semi--direct products of
groups, a construction which translates to crossed products of type $\tto$
factors by groups acting as automorphisms. We also demonstrate that the
hyperfinite type $\tto$ factor can possess both strongly singular masas
and subfactors (see Corollary 2.3, Example 2.5 and also
Section 6). In
the third section we investigate the crossed product of a von Neumann
algebra $\Cal N$ by an action of $\Bbb Z$, with particular
reference to the abelian subalgebra $\Cal B$ generated by the unitary
which implements the action. When the action is ergodic and the
automorphisms are trace preserving, it is well known that the resulting
crossed product is a type $\tto$ factor, [11]. Under the additional
hypothesis that the action is either strongly or weakly mixing (Lemmas
3.1 and 3.2), we obtain strong singularity of $\Cal
B$, showing that such masas arise naturally from classical ergodic theory.

In the remaining three sections we examine the situation of a discrete
I.C.C. group $\Gamma$ with an abelian subgroup $\Gamma_0$, and we consider
$VN(\Gamma_0)$ as an abelian subalgebra of $VN(\Gamma)$. Lemma
4.1 gives a group--theoretic criterion for strong singularity,
which we then verify in various contexts. The unifying theme is to let
$\Gamma$ act on a space $X$ of nonpositive curvature, and to exploit the
geometry to
show that the hypothesis of this lemma is satisfied.
The classical example of $VN({\Bbb{F}}_2)$, alluded to earlier, fits into
this pattern,
by considering the action of ${\Bbb{F}}_2$ on its Cayley graph $X$: a
homogeneous tree of degree four. Now the group ${\Bbb{F}}_2$ is a torsion
free
lattice in the rank one $p$-adic semisimple group
${\hbox{\rm{SL}}}_2({\Bbb{Q}}_p)$.
Conversely any torsion free lattice in ${\hbox{\rm{SL}}}_2({\Bbb{Q}}_p)$ is
a free group
[19].
Any semisimple Lie group acts in a natural way upon a space of nonpositive
curvature.
In the case of a real group this space is a (Riemannian) symmetric space,
whereas in the $p$-adic case it is a euclidean building (a tree, if the group
has rank one).
These cases are dealt with in Sections 4 and 5 respectively.
In Section 4, $\Gamma$ is the fundamental group of a compact locally symmetric
space of nonpositive curvature, while in Section 5, $\Gamma$ acts
by isometries on locally finite euclidean buildings. The motivating
examples here come from [17], which was our starting point in
constructing strongly singular masas.
The building examples are of particular interest because many of the groups
constructed
in [5] do not embed naturally into linear groups.  In the final
section of the paper,
we give another class of examples based on, but extending, those of
Dixmier, [6]. As before, the geometry of the spaces on which our
groups act is the crucial ingredient.

In [16], Popa was able to construct singular masas in any type $\tto$
factor. At this time, we do not know if this is also true for strongly
singular masas, or indeed whether all singular masas must also be strongly
singular.

The third author gratefully acknowledges the hospitality and support
of the Department of Mathematics, University of Newcastle, where part of this
work was completed.

\head 2. Group von Neumann algebras\endhead

Let $\cl A \subseteq \cl M$ be a von Neumann subalgebra of a type $\tto$
factor
$\cl M$. For the case of an abelian subalgebra, {\it {strong singularity}} of
$\cl A$
was defined, in [20], by requiring the inequality
$$
\|\bb E_{\cl A} - \bb E_{u\cl Au^*}\|_{\infty,2} \ge \|(I-\bb E_{\cl A})(u)\|_2
\tag{2.1}
$$
to hold for all unitaries $u\in\cl M$. There is no need for commutativity of
$\cl A$, and so the definition extends without change to all von~Neumann
subalgebras.

For masas $\cl A\subseteq \cl M$, the concept of an asymptotic homomorphism
(with respect to a unitary $v\in\cl A$) was introduced in [20]. We say
that $\bb E_{\cl A}$ is an asymptotic homomorphism, with respect to $v\in\cl
A$,
if
$$
\lim_{|n|\to \infty} \|\bb E_{\cl A}(xv^ny) - \bb E_{\cl A}(x) v^n\bb E_{\cl
A}(y)\|_2 = 0\tag{2.2}
$$
for all $x,y\in\cl M$. Strong singularity is a consequence of having an
expectation which is an asymptotic homomorphism, [20], and
(2.2) gives a criterion for strong singularity which can be easily
checked in specific cases. Our first result is a minor variant of this. We
weaken the requirement of (2.2) slightly, but obtain the same
conclusion. The resulting criterion is then more flexible. The proof is
essentially that given in [20], but we include it for completeness.
It will become apparent later why we state condition (2.3) in a
stronger form
than is necessary for the proof of the result.

\proclaim{Lemma 2.1}
Let $\cl A\subseteq \cl M$ be a von Neumann subalgebra of a type $\tto$ factor
$\cl M$ and suppose that, given $\vp>0$ and $\{x_1,\ldots, x_n; y_1,\ldots,
y_m\}\in\cl M$, there exists a unitary $v\in \cl A$, such that
$$
\|\bb E_{\cl A}(x_iv y_j) - \bb E_{\cl A}(x_i)v \bb E_{\cl A}(y_j)\|_2 < \vp.
\tag{2.3}
$$
Then $\cl A$ is strongly singular in $\cl M$.
\endproclaim

\demo{Proof}
We will make use of the simple relationship
$$
\|h\|^2 = \|Ph\|^2 + \|(I-P)h\|^2\tag{2.4}
$$
for any element $h$ in a Hilbert space $H$ and for any projection $P\in B(H)$.
Fix a unitary $u\in\cl M$ and $\vp>0$. Apply the hypothesis to the set
$\{u^*;u\}$ to obtain a unitary $v\in \cl A$ such that
$$
\|\bb E_{\cl A}(u^*v u) - \bb E_{\cl A}(u^*)v \bb E_{\cl A}(u)\|_2 <
\vp.\tag{2.5}
$$
Using this inequality, we see that
$$\align
\|\bb E_{\cl A} - \bb E_{u\cl Au^*}\|^2_{\infty,2} &\ge \|v-\bb E_{u\cl
Au^*}(v)\|^2_2\\
&= \|v-u\bb E_{\cl A}(u^*vu)u^*\|^2_2\\
&= \|u^*vu-\bb E_{\cl A}(u^*vu)\|^2_2\\
&= 1 - \|\bb E_{\cl A}(u^*vu)\|^2_2\\
&\ge 1 - (\|\bb E_{\cl A}(u^*)v\bb E_{\cl A}(u)\|_2 + \vp)^2\\
&\ge 1-(\|\bb E_{\cl A}(u)\|_2 + \vp)^2 \\
&= \|(I-\bb E_{\cl A})(u)\|^2_2 - \vp^2-2\vp\|\bb E_{\cl A}(u)\|_2.\tag{2.6}
\endalign$$
Since $\vp>0$ was arbitrary, the result follows.
\enddemo

The two basic ways to  obtain type $\tto$ factors are to consider the
von~Neumann algebras arising from discrete groups, and to form crossed
products
by group actions. Such actions on group von~Neumann algebras can take place at
the level of groups, as we now describe. Let $K$ and $H$ be discrete groups
with
an action $\alpha\colon \ K\to \hbox{Aut}(H)$. Then the semi-direct product
${\Gamma}=H
\rtimes_\alpha K$ is the set of formal products $\{hk\colon \ h\in H,\ k \in
K\}$
with multiplication
$$
(hk)(h'k') = (h\alpha_k(h'))(kk').\tag{2.7}
$$
The action $\alpha$ lifts from $H$ to $VN(H)$, and $VN({\Gamma}) = VN(H)
\rtimes_\alpha K$,
[21]. We assume this notation in the next result. Identity elements of
groups are denoted
$e_H$ or $e_K$, and the abbreviation ${\hbox{I.C.C.}}$ means {\it{infinite
conjugacy
class}}.

\proclaim{Theorem 2.2}
Let $H$ and $K$ be infinite discrete groups, let $\alpha\colon \ K\to
\hbox{Aut}(H)$ be an action, and let ${\Gamma} = H\rtimes_\alpha K$. Consider
the following
statements.
\roster
\item"{(i)}"
For each $k\in K\backslash\{e_K\}$, the only fixed point of
$\alpha_k$ is $e_H$;
\item"{(ii)}"
 ${\Gamma}$ is I.C.C.\ and, given finite subsets
$F_1,F_2\subseteq
H\backslash\{e_H\}$, there exists $k\in K$ such that $\alpha_k(F_1)\cap F_2 =
\emptyset$;
\item"{(iii)}"
 ${\Gamma}$ is I.C.C.\ and, given a finite subset
$F\subseteq
H\backslash
\{e_H\}$, there exists $k\in K$ such that $\alpha_k(F) \cap F = \emptyset$;
\item"{(iv)}"
 $VN({\Gamma})$ is a type $\tto$ factor and $VN(K)$ is a strongly
singular
von~Neumann subalgebra.
\endroster
Then (i) $\Rightarrow$ (ii) $\Leftrightarrow$ (iii) $\Rightarrow$ (iv).
\endproclaim

\demo{Proof}
The equivalence of (ii) and (iii) is clear. In one direction, given $F_1$ and
$F_2$, take $F$ to be $F_1\cup F_2$; in the other, given $F$, take $F_1 = F_2=
F$. We now show that (i) implies (iii). Suppose that (i) holds, but that there
exists a finite set $F\subseteq H\backslash \{e_H\}$ so that $\alpha_k(F) \cap
F\ne \emptyset$ for all $k\in K$. Then let
$$
S_f = \{k\in K\colon \ \alpha_k(f)\in F\},\qquad f\in F,\tag{2.8}
$$
and note that $K = \bigcup\limits_{f\in F} S_f$, so that at least one $S_f$ is
infinite. For this $f$, there must exist distinct $k_1,k_2\in K$ such that
$\alpha_{k_1}(f) = \alpha_{k_2}(f)$, since $F$ is finite. But then
$\alpha_{k_1k_2^{-1}}$ has a fixed point, contradicting (i). It remains to
show
that ${\Gamma}$ is I.C.C.

Consider first $k_1\in K\backslash \{e_K\}$. Then
$$
h^{-1}k_1h = h^{-1}\alpha_{k_1}(h)k_1,\qquad h\in H,\tag{2.9}
$$
and so there are infinitely many distinct conjugates of $k_1$ unless
$\alpha_{k_1}$ has fixed points, contrary to hypothesis. If $h_1\ne e_H$,
consider conjugates of $h_1k_1\in {\Gamma}$ by elements $k\in K$. Then
$$
kh_1k_1k^{-1} = \alpha_k(h_1) kk_1k^{-1},\qquad k\in K.\tag{2.10}
$$
The set $\{\alpha_k(h_1)\colon \ k\in K\}$ consists of an infinite number of
distinct elements, otherwise $h_1$ is a fixed point of some $\alpha_k$, and so
${\Gamma}$ is I.C.C.

We now assume (iii). Then the I.C.C. condition on ${\Gamma}$ ensures that
$VN({\Gamma})$ is a
type $\tto$ factor. We will verify that (2.3) is satisfied, and then
obtain the result from Lemma~2.1. A simple approximation argument
shows
that we may take the $x_i$'s and $y_j$'s to be group elements in (2.3).
Moreover, by expanding the set of such elements, we may assume that the
inequality to be verified takes the form
$$
\|\bb E_{\cl A}(x_ivx_j) - \bb E_{\cl A}(x_i)v\bb E_{\cl A}(x_j)\|_2 < \vp
\tag{2.11}
$$
for a given set $\{x_1,\ldots, x_r\} \subseteq {\Gamma}$. The module map
property of $\bb
E_{\cl A}$ shows that (2.11) is true for any $x_i\in K$, so we may
assume
that each $x_i$ has the form $h_ik_i$ with $k_i\in K$ and $h_i\in H\backslash
\{e_H\}$. Then $\bb E_{\cl A}(x_i) = 0$, $1\le i \le r$, so (2.11)
will
be satisfied by a group element $k\in K$, chosen so that
$$
h_ik_ikh_jk_j\notin K,\qquad 1\le i,j\le r.\tag{2.12}
$$
This condition is equivalent to
$$
h_i\alpha_{k_ik}(h_j) k_ikk_j\notin K,\qquad 1\le i,j\le r,\tag{2.13}
$$
which will be true if
$$
h_i\alpha_{k_ik} (h_j) \ne e_H,\qquad 1\le i,j\le r.\tag{2.14}
$$
This last condition may be reformulated as
$$
\alpha_k(h_j)\ne \alpha_{k^{-1}_i}(h^{-1}_i),\qquad 1\le i,j\le r.\tag{2.15}
$$
Let $F = \{h_1,\ldots, h_r, \alpha^{-1}_{k_1}(h^{-1}_1),\ldots,
\alpha_{k^{-1}_r} (h^{-1}_r)\} \subseteq H\backslash \{e_H\}$. By hypothesis,
there exists $k\in K$ such that $\alpha_k(F)\cap F = \emptyset$. In particular
(2.15) is satisfied for this choice of $k$, completing the proof.
\enddemo

It is now easy to produce examples of strongly singular subalgebras by
constructing groups which satisfy Theorem~2.2 (i).

\proclaim{Corollary 2.3}
The hyperfinite type $\tto$ factor contains a strongly singular hyperfinite
subfactor.
\endproclaim

\demo{Proof}
Let $K$ be a countable amenable I.C.C.\ group with no elements of finite order
except the identity. An example of such a group is given below. Then let $H$
be
the countable abelian group, under pointwise multiplication, of functions
$f\colon \ K\to \{\pm 1\}$ which are identically 1 off a finite set. Then
define
an action $\alpha\colon \ K\to \hbox{Aut}(H)$ by
$$
\alpha_k(f)(x) = f(k^{-1}x),\qquad x\in K,\tag{2.16}
$$
for each $k\in K$ and $f\in H$. Consider a fixed $k\in K\backslash \{e_K\}$,
and
suppose that $f\in H\backslash \{e_H\}$ is a fixed point of $\alpha_k$, and is
thus a fixed point for all powers of $\alpha_k$. There exists $x_0\in K$ such
that $f(x_0) = -1$, and it then follows that
$$
f(k^{-n}x_0) = -1,\qquad n\in \bb Z.\tag{2.17}
$$
The definition of $H$ shows that $\{k^{-n}x_0\colon \ n\in \bb Z\}$ is a
finite
set, contradicting the assumption that $k$ has infinite order.

Let ${\Gamma} = H\rtimes_\alpha K$. Then $VN({\Gamma}) = VN(H)\rtimes_\alpha
K$ and so is hyperfinite. The
hypothesis of Theorem~2.2 (i) is satisfied, and so $VN({\Gamma})$ is
the
hyperfinite type $\tto$ factor, while $VN(K)$ is a strongly singular
hyperfinite
type $\tto$ subfactor, by (i) $\Rightarrow$ (iv) of this theorem.
\enddemo

\example{Example 2.4}
Take the group $\bb Z$, and let $H$ be the group, under pointwise addition, of
functions $f\colon \ \bb Z\to\bb Z$ which are identically 0 off a finite set.
Define an action $\alpha\colon \ \bb Z\to \hbox{Aut}(H)$ by
$$
\alpha_k(f)(n) = f(n+k),\qquad n\in \bb Z,\tag{2.18}
$$
for $k\in\bb Z$ and $f\in H$. For $k\ne 0$, $\alpha_k$ has no fixed points
other than the identity of $H$, and
so (i) $\Rightarrow$ (ii) of Theorem~2.2 shows that
$K=H\rtimes_\alpha\bb Z$
is I.C.C. Since $H$ and $\bb Z$ are abelian it follows that the semi-direct
product is amenable. It is easy to verify that elements of finite order ($\ne
e$) in
semi-direct products can exist only when they exist in at least one of the
constituent
subgroups, so the group $K$ defined above is an example of an amenable I.C.C.\
group with no elements of finite order except the identity. Another
possibility is to
let the multiplicative group ${\bb Q}^+$ act on the additive group $\bb Q$
by $\alpha_q(p)=qp$. The resulting semi--direct product has exactly the same
properties.$\hfill\square$
\endexample

\example{Example 2.5}
With the notation of the previous example, Theorem~2.2 shows that
$VN(\bb Z)$ is a strongly singular masa inside the hyperfinite type $\tto$
factor $VN(K)$.$\hfill\square$
\endexample

\example{Example 2.6}
Let $\bb F_\infty$ be the free group on countably many generators $g_i$,
indexed
by $i\in \bb Z$. For each $n\in\bb Z$, the map $\alpha_n\colon \ g_i\mapsto
g_{i+n}$ on generators lifts to an automorphism of $\bb F_\infty$ with no
fixed
points except the identity when $n\ne 0$. Theorem~2.2 (i) then gives
a
strongly singular masa $VN(\bb Z)$ inside the non-hyperfinite factor $VN(\bb
F_\infty) \rtimes_\alpha\bb Z$.$\hfill\square$
\endexample

\head 3. Ergodic actions\endhead

In this section we consider a von Neumann algebra ${\cl N}$ (denoted by ${\cl
A}$ when
abelian) with a faithful normal bounded trace $tr$, together with a trace
preserving automorphism $\theta$. We assume that ${\cl N}$ is represented on
$L^2({\cl N},tr)$, and we define an action of $\bb Z$ on ${\cl N}$ by
$\alpha_n =
\theta^n$. The resulting crossed product ${\cl N}\rtimes_\alpha\bb Z$ is
represented on
$L^2({\cl N},tr)\otimes \ell^2(\bb Z)$. There is a representation $\pi$ of
${\cl N}$ on this
Hilbert space and a unitary operator $u$ so that
$$
\pi(\theta^n(x)) = u^nxu^{-n},\qquad x\in {\cl N},\quad n\in\bb Z,\tag{3.1}
$$
and elements of ${\cl N}\rtimes_\alpha\bb Z$ have unique representations as
$\sum\limits_{n\in\bb Z}\pi(x_n)u^n$, where such sums converge ultraweakly.
Since $\theta$ is trace preserving, there is a faithful normal trace on
${\cl N}\rtimes_\alpha \bb Z$ given by
$$
tr\left(\sum_{n\in\bb Z} \pi(x_n)u^n\right) = tr(x_0).\tag{3.2}
$$
This is standard theory which may be found in [11]. When ${\cl N}$ is an
abelian
von~Neumann algebra ${\cl A}$, its image in ${\cl A}\rtimes_\alpha\bb Z$ has a
normalizer which
generates the crossed product. Thus ${\cl A}$ is Cartan whenever the action is
such
that ${\cl A}$ is maximal abelian in ${\cl A}\rtimes_\alpha\bb Z$. The unitary
$u$ which
implements $\theta$ always generates a canonical abelian von~Neumann
subalgbebra
${\cl B}$ of ${\cl A}\rtimes_\alpha\bb Z$, which we now investigate.

Recall that an action $\alpha$ of $\bb Z$ on ${\cl N}$ is {\it {ergodic}} if
each
$\alpha_n$ ($n\ne 0$) has only multiples of the identity as its fixed points.
If
the automorphism group $\{\theta^n\}_{n\in\bb Z}$ is both trace preserving and
ergodic then it is well known (see [11, p.~546]) that ${\cl
N}\rtimes_\alpha\bb Z$ is a
type $\tto$ factor with ${\cl B}$ as a masa. The automorphism $\theta$ is {\it
{strongly
mixing}} (called {\it {mixing}} in classical abelian ergodic theory) if
$$
\lim_{n\to\infty} tr(x\theta^n(y)) = tr(x) tr(y),\qquad x,y\in {\cl
N}.\tag{3.3}
$$
When $\theta$ is trace preserving, we may use the limit as $|n|\to\infty$ in
(3.3), and we also note that ergodicity is an easy consequence of
(3.3). We say that $\theta$ is {\it {weakly mixing}} if
$$
\lim_{n\to\infty} \frac1n \sum^{n-1}_{k=0} |tr(x\theta^k(y)) - tr(x)tr(y)| = 0
\tag{3.4}
$$
for all $x,y\in {\cl A}$. Ergodicity also follows from this weaker definition.

\proclaim{Lemma 3.1}
Let ${\cl N}$ be a von Neumann algebra with a faithful normal trace, let
$\theta$ be a
trace preserving ergodic automorphism of ${\cl N}$ defining an action
$\alpha_n =
\theta^n$ of $\bb Z$ on ${\cl N}$, let $u$ be the unitary in ${\cl M}={\cl
N}\rtimes_\alpha\bb Z$ which
implements $\theta$, and let ${\cl B}$ be the abelian von~Neumann algebra
generated by
$u$. Then $\bb E_{\cl B}$ is an asymptotic homomorphism with respect to $u$ if
and
only if $\theta$ is strongly mixing. In particular, ${\cl B}$ is a strongly
singular
masa in ${\cl M}$ when $\theta$ is strongly mixing.
\endproclaim

\demo{Proof}
The set $\{u^n\}_{n\in\bb Z}$ is an orthonormal basis for $L^2({\cl B},tr)$,
and so
$\bb E_{\cl B}$ may be expressed by
$$
\bb E_{\cl B}(x) = \sum_{n\in\bb Z} tr(xu^{-n})u^n,\qquad x\in {\cl
M}.\tag{3.5}
$$
In particular
$$
\bb E_{\cl B}(x) = tr(x)1,\qquad x\in {\cl N}.\tag{3.6}
$$
Elements  of the form $xu^n$, $x\in {\cl N}$, $n\in\bb Z$, generate ${\cl M}$,
so it is
sufficient to check the asymptotic homomorphism condition for such operators.
If
$x,y\in {\cl N}$ and $k,r\in\bb Z$, then
$$\align
\bb E_{\cl B}(xu^ku^nyu^r) &= \bb E_{\cl
B}(x\theta^{n+k}(y)u^{n+k+r})\\
&= \bb E_{\cl B}(x\theta^{n+k}(y)) u^{n+k+r}\\
&=  tr(x\theta^{n+k}(y))u^{n+k+r},\tag{3.7}
\endalign$$
from (3.6). On the other hand,
$$
\bb E_{\cl B}(xu^k)u^n\bb E_{\cl B}(yu^r) = tr(x) tr(y)u^{n+k+r},\tag{3.8}
$$
using (3.5), and so
$$\align
\|\bb E_{\cl B}(xu^ku^nyu^r) - \bb E_{\cl B}(xu^k)u^n\bb E_{\cl B}(yu^r)\|_2
&=
|tr(x\theta^n(\theta^k(y))) - tr(x)tr(y)|\\
&= |tr(x\theta^n(\theta^k(y))) - tr(x) tr(\theta^k(y))|,\tag{3.9}
\endalign$$
where the last equality uses trace preservation by $\theta$. Letting $|n|\to
\infty$ in (3.9) immediately gives the conclusion that $\bb E_{\cl B}$
is an
asymptotic homomorphism for $u$ if and only if $\theta$ is strongly mixing.
The
last statement of the lemma then follows from [20].
\enddemo

We now consider a weakly mixing automorphism $\theta$, and we maintain the
notation of the previous lemma.

\proclaim{Lemma 3.2}
Let $\theta$ be a trace preserving weakly mixing automorphism of ${\cl N}$.
Then ${\cl B}$
is a strongly singular masa in ${\cl M} = {\cl N}\rtimes_\alpha\bb Z$.
\endproclaim

\demo{Proof}
We will verify (2.3) in Lemma 2.1, from which the result will
follow. It suffices to consider a finite set of generators so, to obtain a
contradiction, we may assume that there exist $\vp>0$ and elements $x_ju^k\in
{\cl M}$, $0\le j$, $|k|\le J$, so that (2.3) fails for all unitaries
$v\in
{\cl B}$. In particular
$$
\max\{\|\bb E_{\cl B}(x_ju^{n+k}x_ru^s) - \bb E_{\cl B}(x_ju^k)u^n\bb E_{\cl
B}(x_ru^s)\|_2 \colon
\ 0 \le j,r,|k|, |s| \le J\} \ge \vp\tag{3.10}
$$
for all $n\in\bb Z$. Using (3.7), this condition becomes
$$
\max\{|tr(x_j\theta^{n+k}(x_r)) - tr(x_j)tr(\theta^k(x_r))|\colon \  0 \le
j,r,
|k|, \le J\} \ge\vp\tag{3.11}
$$
for all $n\in \bb Z$. Let $y_{k,r}$ denote the element $\theta^k(x_r)$. If we
sum in (3.11) and average from $0$ to $n-1$, then we obtain
$$
\frac1n \sum^{n-1}_{i=0} \sum_{0\le j,r,|k|\le J} |tr(x_j\theta^i(y_{k,r}) -
tr(x_j)tr(y_{k,r})| \ge \vp\tag{3.12}
$$
for all $n\ge 1$, and this violates the defining inequality (3.4) of
weakly mixing. This completes the proof .
\enddemo

\remark{Remark \rom{3.3}}
Classical ergodic theory (see [15]) provides many examples of strongly
mixing transformations of measure spaces, as well as examples which are weakly
but not strongly mixing. The two previous lemmas then give examples of
strongly
singular masas ${\cl B}$, some of which do not arise from asymptotic
homomorphisms for
$u$.$\hfill\square$
\endremark

\head 4. Groups acting on symmetric spaces\endhead

Let $\Gamma$ be an I.C.C.\ group with an abelian subgroup
$\Gamma_0$. Then
$VN(\Gamma_0)$ is an
abelian subalgebra of the type $\tto$ factor $VN(\Gamma)$, and in this section
we investigate when
it is a strongly singular masa. For the case of group von~Neumann algebras,
Lemma 2.1
takes the following form.

\proclaim{Lemma 4.1}
Let $\Gamma$ be a discrete I.C.C.\ group with an abelian subgroup
$\Gamma_0$.
The following condition implies that $\VN(\Gamma_0)$ is a strongly
singular masa of $\VN(\Gamma)$: \\
If $x_1,\ldots,x_m \in \Gamma$ and
$$
\Gamma_0\subseteq\bigcup_{i,j}x_i\Gamma_0 x_j\, ,\tag{4.1}
$$
then $x_i\in \Gamma_0$ for some $i$.
\endproclaim

\demo{Proof}
The condition in question is equivalent to the following: \\
If
$x_1,\ldots ,x_m,y_1,\ldots ,y_n \in \Gamma \backslash \Gamma_0$, then there
exists $\gamma_0\in \Gamma_0$ such that
$$
x_i\gamma_0y_j\notin \Gamma_0,\ \ 1 \leq i \leq m,\ 1 \leq j \leq n.\tag{4.2}
$$
To see this replace each of the sets $\{x_1,\ldots, x_n\}$, $\{y_1,\ldots,
y_m\}$ by their union
and replace $x_i$ by
$x_i^{-1}$, and $y_j$ by $y_j^{-1}$.
Now apply Lemma 2.1, with each operator approximated by a finite
linear combination of
group elements.
\enddemo

The aim now is to apply this lemma to construct
strongly singular masas of $VN(\Gamma)$, for
certain geometrically defined groups $\Gamma$, acting on spaces of nonpositive
curvature.

In order to establish a connection with geometry, consider the following
general setup.
Let $(X,d)$ be a metric space and let $\Gamma$ be a group of isometries of
$X$.
If $P$, $Q$ are subsets of $X$, and $\delta>0$, then use the notation
$P\displaystyle\operatornamewithlimits{\subset}_{\delta}Q$ to mean that
$d(p,Q)\leq\delta$,
for all $p\in
P$.

Let $\Gamma_0$ be an abelian subgroup of $\Gamma$ and let $A$ be a
$\Gamma_0$-invariant subset of $X$.  Consider the conditions:
\roster
\item"{\bf (C1)}"
 There exists a compact subset $K$ of $A$ such that $\Gamma_0
K=A$.\\
\item"{\bf (C2)}"
 If $A\displaystyle\operatornamewithlimits{\subset}_{\delta} x_1A\cup
x_2A\cup\cdots\cup
x_mA$,
for some $x_1,\ldots,x_m\in\Gamma$, and $\delta>0$, then $x_j\in\Gamma_0$, for
some $j$.
\endroster

\proclaim{Proposition 4.2}
If {\rm {\bf (C1)}} and {\rm {\bf (C2)}} hold then
$VN(\Gamma_0)$ is a strongly singular masa of $VN(\Gamma)$.
\endproclaim

\demo{Proof}
Suppose that $x_1,\ldots,x_m \in\Gamma$ and
$$
\Gamma_0\subseteq\bigcup_{i,j}x_i\Gamma_0 x_j\,.\tag{4.3}
$$
Let $\delta=\max\{d(x_jk,k);1\leq j\leq n,\,k\in K\}$.
For $1\leq j\leq n$, this implies that $x_j
K\displaystyle\operatornamewithlimits{\subset}_{\delta} K$ and
so
$$
\Gamma_0x_j K\displaystyle\operatornamewithlimits{\subset}_{\delta}\Gamma_0
K=A\,.\tag{4.4}
$$
Hence, for each $i$, $j$, we have
$x_i\Gamma_0x_jK\displaystyle\operatornamewithlimits{\subset}_{\delta}x_iA$.
It follows from {\bf (C1)} and (4.3) that
$$
A=\Gamma_0K\displaystyle\operatornamewithlimits{\subset}_{\delta}x_1A\cup
x_2A\cup\cdots\cup
x_mA\,.\tag{4.5}
$$
Applying condition {\bf (C2)}, we see that $x_j\in\Gamma_0$ for some $j$,
contrary to hypothesis.
The result now follows from Lemma 4.1.
\enddemo

In the first class of examples, $\Gamma$ is the fundamental group of a compact
locally symmetric space of nonpositive curvature.  The classic
book [13] is a convenient reference for the background and
necessary results. There is a clear introduction to the theory of symmetric
spaces in
[2, Chapter II.10].

Let $X$ be a symmetric space of noncompact type, by which we mean
that $X$ is a quotient
$G/K$ of a semisimple Lie group by a maximal compact subgroup $K$.

\proclaim{Lemma 4.3}
$X=\SL_n(\RR)/\SO_n(\RR)$, $n\geq2$.
If $n=2$ one obtains the hyperbolic plane.$\hfill\square$
\endproclaim

The {\it rank} $r$ of $X$ is the dimension of a maximal {\it flat} in $X$.
That is, the
maximal dimension of an isometrically embedded
euclidean space in $X$.  If $r=1$, then the flats are geodesics.
If $X=\SL_n(\RR)/\SO_n(\RR)$, then $r=n-1$.  Call a flat of maximal dimension
$r$
an {\it $r$-flat} (or {\it maximal flat}).
A geodesic $L$ in $X$ is called {\it regular} if it lies in only one
$r$-flat; it is called {\it singular} if it is not regular.

Let $F$ be an $r$-flat in $X$ and let $x\in F$.  Let $S_x$ denote the
union of all the singular geodesics through $x$.  A connected component
of $F-S_x$ is called a {\it Weyl chamber} with origin $x$.

\example{Example 4.4}
If $X=\SL_2(\RR)/\SO_2(\RR)$, the hyperbolic plane, then $r=1$.
The 1-flats are geodesics.  If $x$ is a point on a geodesic $L$
then the two Weyl chambers in $L$ with origin $x$ are two
semi-geodesics.$\hfill\square$
\endexample

\example{Example 4.5}
If $X=\SL_3(\RR)/\SO_3(\RR)$, then $r=2$ and there are six Weyl chambers in
any 2-flat $F$ with a given origin $x\in F$, as illustrated in
Figure 1.
\bigskip

\centerline{
\beginpicture
\setcoordinatesystem units <1cm, 1.732cm>
\setplotarea  x from -2 to 2,  y from -1 to 1
\putrule from -2 0 to 2 0
\setlinear \plot  -1 -1  0 0  1 1  /
\setlinear \plot  1 -1  0 0  -1 1 /
\put {$_{\bullet}$} at 0 0
\put {$x$} [b,l] at 0.3 0.1
\endpicture}\bigskip
\hfil
\centerline{Figure 1. Weyl chambers with origin $x$}

\endexample

If $A$, $B$ are subsets of $X$, define the Hausdorff distance between them by
$$
\hd(A,B)=\inf\{\delta\leq\infty;
A\displaystyle\operatornamewithlimits{\subset}_{\delta}B\hbox{ and
}B\displaystyle\operatornamewithlimits{\subset}_{\delta} A\}\,.\tag{4.6}
$$
Let ${\Cal W}$ denote the set of all Weyl chambers in $X$ and
define an equivalence relation $\sim$ on ${\Cal W}$ by
$$
W_1\sim W_2\Longleftrightarrow\hd(W_1,W_2)<\infty\,.\tag{4.7}
$$
The {\it boundary} $\Omega$ of $X$ is defined to be the quotient space
${\Cal W}/\sim$.  It is well known [13, Lemma 4.1] that
$\Omega$ may be identified with the topological homogeneous space
$G/P$, where $P$ is a Borel subgroup of $G$.  The action
of a discrete subgroup $\Gamma$ of $G$ on the boundary $\Omega$ will
play an important role in our argument, just as it did in Mostow's proof
of rigidity [13].

If $F$ is an $r$-flat in $X$, then the restriction of the equivalence
relation $\sim$ to $F$ allows one to define the boundary of $F$, which is a
finite set. There is a natural embedding of the boundary of $F$ into the
boundary of $X$ and it is convenient to identify each boundary point of $F$
with the corresponding boundary point of $X$.


Suppose that $\Gamma$ is a cocompact lattice in a
semisimple Lie group $G$. It is well known
[2, Proposition II.6.10], [13, \S 11] that
each element $\gamma\in \Gamma$ is {\it semi-simple}.
Geometrically this means
that the displacement function defined on $X$ by $\xi \mapsto d(\xi,
\gamma\xi)$
attains its minimum at some point $\xi_0\in X$. If the group $\Gamma$ acts
freely
on $X$ then the minimum value $d(\xi_0, \gamma\xi_0)$ is strictly positive if
$\gamma\ne 1$ ($\gamma$ is {\it hyperbolic}). This implies [2, Proposition
II.6.8]
that there is a geodesic line (an {\it axis} of $\gamma$)
upon which $\gamma$ acts by translation.


We now have enough background information to begin the
main result of this section. Let $G$ be a semisimple Lie group with no centre
and no compact factors.
Let $\Gamma$ be a torsion free cocompact lattice in $G$.  Then $\Gamma$  acts
freely
on the symmetric space $X=G/K$ and the quotient manifold $M=\Gamma\backslash
X$
has universal covering space $X$.  Thus $M$ is a compact locally symmetric
space of nonpositive curvature, with fundamental group $\pi (M)=\Gamma$.
Moreover every compact locally symmetric space $M$ arises in this way.

Let $T^r\subset M$ be a totally geodesic embedding of a flat $r$-torus
in $M$.  By the easy part of the Flat Torus Theorem [12, Theorem 1]
the inclusion $i:T^r\to M$ induces an injective homomorphism
$i_*:\pi(T^r)\to\pi(M)$.
Thus $\Gamma_0=i_*\pi(T^r)\cong\ZZ^r$.
Conversely, if $\Gamma_0$ is any free abelian subgroup of rank $r$ in
$\Gamma$,
then by [12, Theorem 1], [2, Theorem II.7.1], there exists
an $r$-flat $F_0$ in $X$ such that $\Gamma_0F_0=F_0$, $\Gamma_0$ acts
upon $F_0$ by translations, and $\Gamma_0\backslash F_0=T^r$.

Let $\sigma=\sigma(M)$ denote the length of a shortest closed geodesic
in $M$.  The aim of the rest of this section is to prove strong singularity of
$VN(\Gamma_0)$
in this setting. We accomplish this through the following series of lemmas.
Note that the group $\Gamma$ is ${\hbox{I.C.C.}}$ by [9, Lemma 3.3.1],
so
$\VN(\Gamma)$
is a $\tto$ factor.
We shall verify conditions {\bf (C1)}, {\bf (C2)} for the action of
$\Gamma$ on the symmetric space $X$, taking the subset $A$ of $X$ to be
the $r$-flat $F_0\subset X$, upon which the abelian subgroup
$\Gamma_0\cong\ZZ^r$ acts.
The result is then a consequence of Proposition 4.2.
Verification of {\bf (C1)} is easy.

\proclaim{Lemma 4.6}
The action of $\Gamma_0$ on $F_0$ satisfies {\rm{\bf (C1)}}.
\endproclaim

\demo{Proof}
This is immediate since $\Gamma_0\backslash F_0$ is compact.  Let $K$ be the
closure of a
bounded fundamental domain for the action of $\Gamma_0$ on $F_0$.
\enddemo

Verification of {\bf (C2)} requires some preparation.

\proclaim{Lemma 4.7}
If $F,F_1,\ldots,F_m$ are $r$-flats in $X$ and
$F\displaystyle\operatornamewithlimits{\subset}_{k}F_1\cup\cdots\cup F_m$,
for some $k>0$, then each boundary point of $F$ is a boundary point of some
$F_j$, $1\leq j\leq m$.
\endproclaim

\demo{Proof}
Let $W$ be a Weyl chamber in $F$.  Write each flat $F_j$\, $(1\leq j\leq m)$,
as a finite union of Weyl chambers $W_{jl}$.  Then $W\subset
\bigcup_{j,l}W_{jl}$.
It follows from [13, Lemma 15.1] that $\hd(W,W_{rs})<\infty$ for some
$r$, $s$.
That is, $W$ and $W_{rs}$ represent the same boundary point of $X$.
\enddemo

\proclaim{Lemma 4.8}
Assume that $T^r$
has diameter $<\sigma$.
If $x\in\Gamma$ and the $r$-flats $F_0$ and $xF_0$ have a common boundary
point,
then $x\in\Gamma_0$.
\endproclaim

\demo{Proof}
This depends crucially upon the fact that the embedded torus
$T^r=\Gamma_0\backslash F_0$
has diameter $<\sigma$, where $\sigma$ is the minimum length of a nontrivial
closed geodesic in $X$.  This means that for any two points $a$, $b\in F_0$,
there exists $\gamma\in\Gamma_0$ such that
$$
d(a,\gamma b)\leq\diam(T^r)<\sigma\,,\tag{4.8}
$$
where $d$ is the canonical $G$-invariant metric on $X$.

The hypothesis that the $r$-flats $F_0$, $xF_0$ have a common boundary point
implies that there exist Weyl chambers $W$, $W'$ in $F_0$, $xF_0$ respectively
such that $\hd(W,W')<\infty$.  Let $v$ be the origin of $W$.  The Weyl chamber
$x^{-1}W'$ in $F_0$ is equivalent to a Weyl chamber $W_0$ in $F_0$ with origin
$v$.
Thus the Weyl chamber $xW_0$ in $xF_0$ is equivalent to $W'$ and hence to $W$.

Choose $\gamma_0\in\Gamma_0$ such that $\gamma_0v\in W_0$.  This is possible,
since $\Gamma_0$ acts freely on $F_0$ by translations and $\Gamma_0\backslash
F_0 =T^r$.
Let $[v,\gamma_0v]$ denote the geodesic segment in $F_0$ from $v$ to
$\gamma_0v$.  Then
$L=\bigcup_{n\geq 0}\gamma_0^n[v,\gamma_0v]$ is a regular geodesic ray
contained in
$W_0$ and $\gamma_0$ acts on $L$ by translation.

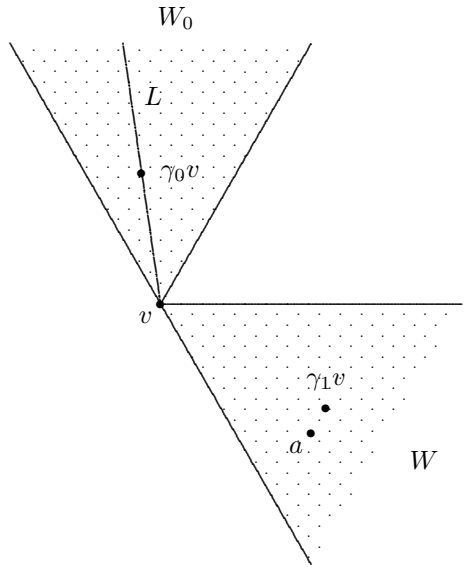
\centerline{
\beginpicture
\setcoordinatesystem units <1cm, 1.732cm>
\setplotarea  x from -2 to 2,  y from -2 to 3
\setlinear \plot  -2 2  0 0  2 2 /
\setlinear \plot  0 0  -0.5 2 /
\put {$_{\bullet}$} at 0 0
\put {$v$}  at -0.2 -0.1
\put {$L$}  at -0.1 1.6
\put {$W_0$}  at 0.2 2.2
\put {$\gamma_0v$} [l] at 0 1
\put {$_{\bullet}$} at -0.25 1
\setshadegrid span <4pt>
\vshade  -2 2 2 <,,,z>  0 0 2  <z,,,>  2 2 2 /
\setlinear \plot  2 -2  0 0  4 0 /
\setshadegrid span <4pt>
\vshade  0 0 0 <,,,z> 2 -2 0  <z,,,>  4 0 0 /
\put {$_{\bullet}$} at 2.2 -0.8
\put {$\gamma_1v$}  at 2.2 -0.6
\put {$_{\bullet}$} at 2 -1
\put {$a$}  at  1.8 -1.1
\put {$W$}  at 3.5 -1.2
\endpicture
}\bigskip
\hfil
\centerline{Figure 2. Weyl chambers with origin $v$ in the flat $F_0$}
\bigskip

It follows that $xL$ is a regular geodesic ray lying in the Weyl chamber
$xW_0$, and
$x\gamma_0x^{-1}$ acts by translation on $xL$.  Thus
$xL=\bigcup_{n\geq 0}[x\gamma_0^nv,x\gamma_0^{n+1}v]$.

Since $\hd(xW_0,W)<\infty$, we have
$xL\displaystyle\operatornamewithlimits{\subset}_{t}W$
for some $t>0$.

It follows that $xL$ is actually asymptotic to $W$, meaning that
 $d(xL,W)=0$ [13, Lemma 7.3(iii)].
Now $xL$ approaches $W$ monotonically at infinity [13, Lemma 4.2] and so
$d(x\gamma_0^nv,W)\to 0$ as $n\to\infty$.

Choose $k\geq 1$ and $a\in W$ such that $d(x\gamma_0^kv,a)<\sigma-\diam(T^r)$.

Using (4.8), choose $\gamma_1\in\Gamma_0$ such that
$d(a,\gamma_1v)\leq\diam(T^r)$.
Then $d(x\gamma_0^kv,\gamma_1v)<\sigma$.  Equivalently,
$d(\gamma_1^{-1}x\gamma_0^kv,v)<\sigma$.

Now this implies that $\gamma_1^{-1}x\gamma_0^kv=v$.
For otherwise the geodesic segment from $v$ to $\gamma_1^{-1}x\gamma_0^kv$ in
$X$ projects to a nontrivial closed geodesic in $M$ of length $<\sigma$,
contradicting the definition of $\sigma$.  Since $\Gamma$ acts freely on $X$,
we deduce that $\gamma_1^{-1}x\gamma_0^k=1$.  Therefore
$x=\gamma_1\gamma_0^{-k}\in\Gamma_0$.
\enddemo

\proclaim{Theorem 4.9}
Let $T^r$ be a totally geodesic flat torus in a compact locally
symmetric space $M$ of nonpositive curvature and rank $r$.
Let $\Gamma_0\cong\ZZ^r$ be the image of the fundamental group
$\pi(T^r)$ under the natural monomorphism from  $\pi(T^r)$ into
$\Gamma=\pi(M)$.
Assume that $\diam(T^r)<\sigma(M)$.  Then $\VN(\Gamma_0)$ is a
strongly singular masa of $\VN(\Gamma)$.
\endproclaim

\demo{Proof}
Condition {\bf (C1)} is satisfied by Lemma 4.6.
It remains to verify condition {\bf (C2)}.  Suppose therefore that
$x_1,\ldots,x_m\in\Gamma$ and $\delta>0$ satisfy
$F_0\displaystyle\operatornamewithlimits{\subset}_{\delta}x_1F_0\cup
x_2F_0\cup\cdots\cup
x_mF_0$.
Choose a boundary point $\omega$ of $F_0$.  By Lemma 4.7, $\omega$ is
also a
boundary point of $x_jF$ for some $j$.  It follows from Lemma 4.8
that $x_j\in\Gamma_0$.
Therefore condition {\bf (C2)} is satisfied.

Finally, $\VN(\Gamma_0)$ is a strongly singular masa of $\VN(\Gamma)$,
by Proposition 4.2.
\enddemo

If $\Gamma$ is a torsion free cocompact lattice in $\psl_2(\RR)$ then the
result of Theorem 4.9
 becomes
particularly simple.

\proclaim{Corollary 4.10}
Let $\Gamma$ be the fundamental group of a compact Riemann surface $M$ of
genus $g\geq 2$.
Let $\gamma_0\in\Gamma$ be the class of a closed geodesic of minimal length in
$M$, and
let $\Gamma_0\cong\ZZ$ be the subgroup of $\Gamma$ generated by $\gamma_0$.
Then $\VN(\Gamma_0)$ is a strongly singular masa of $\VN(\Gamma)$.
\endproclaim

\demo{Proof}
If $C$ is a closed geodesic of minimal length $\sigma$ in the class of
$\gamma_0$,
then $\diam(C)=\frac{1}{2}\sigma<\sigma$.  The result follows directly from
Theorem 4.9.
\enddemo

\remark{Remark \rom{4.11}}
The usual presentation of the fundamental group of the compact Riemann surface
$M$ is as the one-relator group
$$\Gamma=\left\langle
a_1,\ldots,a_g,b_1,\ldots,b_g\left|\prod_{i=1}^g[a_i,b_i]=1\right.\right\rangle
$$
where $[a_i,b_i]=a_ib_ia_i^{-1}b_i^{-1}$.

With this presentation of $\Gamma$, the generator $\gamma_0$ of $\Gamma_0$ can
be any one of the generators $a_i^{\pm 1}$, $b_j^{\pm 1}$.  To see this recall
that, by a theorem of Poincar\'e [10, VB], there exist hyperbolic
isometries $a_i$, $b_j\in\psl_2(\RR)$, $(1\leq i,j\leq g)$, which generate
$\Gamma$ inside $\psl_2(\RR)$.  Moreover one can ensure that a fundamental
domain for the action of $\Gamma$ on the hyperbolic plane $X$ is a regular
hyperbolic $(4g)$-gon $P_g$, and the isometries $a_i$, $b_j$ map $2g$ of the
edges of $P_g$ to the other $2g$ edges in an appropriate way.  To see that we
can choose $\gamma_0=a_1$, for example, choose $v$ to be the mid-point of the
edge of $P_g$ such that $a_1v$ is also the mid-point of an edge (Figure
 3).  Then the geodesic segment $[v,a_1v]$ in $P_g$ projects to a
closed geodesic of minimal length in $M=\Gamma \backslash X$.
\vfill\eject

\centerline{
\beginpicture
\setcoordinatesystem units  <0.5in, 0.25in>   
\setplotarea x from -6 to 6, y from -2 to 15  
\putrule from  0 5.5 to 0 3.5
\arrow <10pt> [.2, .67] from  0 3.5 to  0 3.2
\circulararc 140 degrees from -1.3 0.2 center at  -1 0.4
\circulararc 140 degrees from 0.7 0.2 center at  1 0.4
\circulararc 72 degrees from -0.8 -0.05 center at  -1 -0.6
\circulararc 72 degrees from 1.2 -0.05 center at  1 -0.6
\setquadratic
\plot
2 9   1.8 10.5   2 12 /
\plot
2 12   1.4  12.8   1 14 /
\plot
1  14   0 13.6    -1 14 /
\plot
-1 14   -1.4  12.8  -2 12 /
\plot
-2 12  -1.8  10.5  -2 9  /
\plot
-2 9  -1.4 8.2  -1 7 /
\plot
-1 7   0 7.4  1 7  /
\plot
1 7  1.4  8.2  2 9  /
\plot
-2 0  -1.95 0.4  -1.9 0.55  -1.8 0.8   -1.5 1.2  -1 1.5  -0.5 1.4  -0.2 1.1
0 1
0.2 1.1   0.5 1.4   1 1.5  1.5 1.2   1.8 0.8   1.9 0.55  1.95 0.4  2 0  /
\plot
-2 0  -1.95 -0.4  -1.9 -0.55  -1.8 -0.8   -1.5 -1.2  -1 -1.5  -0.5 -1.4  -0.2
-1.1
0 -1
0.2 -1.1   0.5 -1.4   1 -1.5  1.5 -1.2   1.8 -0.8   1.9 -0.55  1.95 -0.4  2 0 /
\setplotsymbol({$\cdot$}) \plotsymbolspacing=4pt
\setquadratic
\plot
1.8  10.5   0.4 11  0 13.6  /
\setquadratic
\plot
1.2 0.2   1.6 0.6    2  0.1  /
\put{$M$}  at    -3 0
\put{$v$}[l]  at    2 10.5
\put{$a_1v$} [b]  at    0 13.8
\put{$P_2$}  at    -3 10.5
\endpicture
}
\hfil
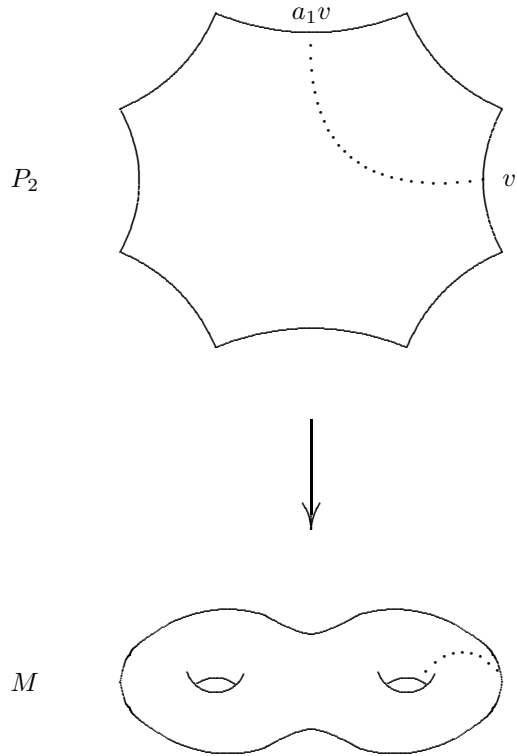
\centerline{Figure 3. The genus 2 case.}

\endremark

\remark{Remark \rom{4.12}}
Corollary 4.10 also follows immediately from [20, Corollary
6.3],
since $\gamma_0$ is a prime element of the non-elementary hyperbolic group
$\Gamma$.$\hfill\square$
\endremark

\head 5. Groups acting on euclidean buildings\endhead

In a second class of examples, the group $\Gamma$ acts cocompactly by
isometries on a
locally finite euclidean building $\Delta$ of rank $r$.  The building $\Delta$
is the
combinatorial counterpart of a symmetric space $X$.  The analogy becomes
particularly
evident if one considers groups of $p$-adic type.  Specifically, let $G$ be a
connected
semisimple group defined over a nonarchimedean local field.  Then $G$ acts  on
its
Bruhat-Tits building $\Delta$ [4], and the vertex set of $\Delta$ may be
identified
with $G/K$, where $K$ is a maximal compact subgroup.

We refer to [18] for the general theory of buildings.
It is worth making a few remarks about the structure of euclidean buildings.

A building $\Delta$ is an $r$-dimensional simplicial complex whose
maximal simplices are called {\it chambers}.  All chambers have the
same dimension $r$ and adjacent chambers have a common face of
dimension $r-1$.

\centerline{
\beginpicture
\setcoordinatesystem units  <1cm, 1.732cm>        
\setplotarea x from -3 to 3, y from -1 to 1.5    
\putrule from 0 0  to  2 0
\setlinear
\plot  2 0   1 1   0 0 /
\plot  0 0   0.8 -0.8   2 0  1.2  -0.8   0 0 /
\plot  0 0   -1 0.8   1 1   -0.8 1.2   0 0  /
\plot  2 0   3 0.8   1 1   2.8 1.2   2 0  /
\endpicture
}\hfil
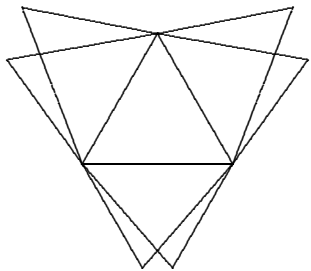
\centerline{Figure 4. Six chambers adjacent to a chamber in an $\tilde A_2$
building.}
\bigskip

 Any two chambers can be connected by a sequence
of adjacent chambers (called a {\it gallery}).  An {\it apartment}
in $\Delta$ is a subcomplex which is isomorphic to a Coxeter complex.
All the apartments are isomorphic and any two simplices in $\Delta$
lie in a common apartment.  If the apartments are infinite then
$\Delta$ is contractible as a topological space.  The apartments
are then euclidean Coxeter complexes isometric to $\RR^n$ and
$\Delta$ is said to be a {\it euclidean building}.
A euclidean building has a canonical piecewise smooth metric which
is consistent with the euclidean structure on the
apartments [3, VI.3]. It is convenient to normalise the
distance on $\Delta$ so that any point of $\Delta$ is at
distance $<1$ from some vertex.
The simplest examples of euclidean buildings are the homogeneous trees.
In such a tree, a chamber is an edge and an apartment is an infinite geodesic.

The boundary of $\Delta$ is defined in terms of equivalence classes of
sectors [18, Chp. 9.3].  A {\it sector} is a simplicial cone of
dimension $r$, with a  {\it special} base vertex, lying in some
apartment of $\Delta$.  A sector in a euclidean building plays the role
of a Weyl chamber in a symmetric space.

\bigskip

\centerline{
\beginpicture
\setcoordinatesystem units  <0.5cm, 0.866cm>        
\setplotarea x from -2.5 to 5, y from -1 to 2.5   
\putrule from -2.5   2     to  2.5  2
\putrule from  -3.5 1  to 3.5 1
\putrule from -4.5 0  to 4.5 0
\putrule from -3.5 -1  to  3.5 -1
\setlinear
\plot  -4.3 -0.3  -1.7 2.3 /
\plot -3.3 -1.3  0.3 2.3 /
\plot -1.3 -1.3  2.3 2.3 /
\plot  0.7 -1.3  3.3 1.3 /
\plot  2.7 -1.3  4.3 0.3 /
\plot  -4.3 0.3  -2.7 -1.3 /
\plot  -3.3 1.3   -0.7 -1.3 /
\plot  -2.3 2.3   1.3 -1.3 /
\plot  -0.3 2.3   3.3 -1.3 /
\plot  1.7 2.3   4.3 -0.3 /
\setshadegrid span <1.5pt>
\vshade  0  0 0 <,z,,>  2.3  0 2.3  <z,,,>  4.3  0 0 /
\endpicture
}\hfil\bigskip
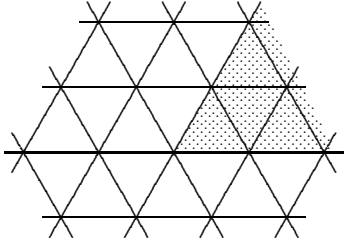
\centerline{Figure 5. Part of an apartment and a sector in an $\tilde A_2$
building.}
\bigskip

Two sectors are said to be
{\it equivalent} if the Hausdorff distance between them is finite.  This is
considerably stronger than the notion of equivalence of Weyl chambers in a
symmetric space, since two sectors are equivalent if and only if they
contain a common subsector.  Equivalent Weyl chambers, by contrast, usually
have no points in common. As for symmetric spaces, the {\it boundary} $\Omega$
of $\Delta$ is the quotient space, whose points are equivalence classes of
sectors.


For the rest of this section we fix a group $\Gamma$ of automorphisms of
$\Delta$ with the following properties.

\roster
\item"{\bf(B1)}" $\Gamma$ acts freely on the vertex set $\Delta^0$, with
finitely many vertex orbits
(i.e. cocompactly).
\item"{\bf (B2)}" There is an apartment $F_0$ in $\Delta$ and an abelian
subgroup $\Gamma_0$ of $\Gamma$ such that $\Gamma_0\bs F_0^0$ is finite,
where $F_0^0$ is the vertex set of $F_0$.\\
\item"{\bf (B3)}"  The natural mapping $\Gamma_0\bs F_0^0\to \Gamma\bs
\Delta^0$ is injective.
\endroster

\remark{Remark \rom{5.1}}
\roster
\item"{(a)}" These conditions have been stated in a form applicable to a
large class of examples.  The assumption {\bf (B1)} does not imply that
$\Gamma$ acts freely on (the geometric realization of) $\Delta$.  In fact,
$\Gamma$ acts freely on $\Delta$ if and only if $\Gamma$ is torsion free
[8, Proof of Theorem 4.1].

If we assume that $\Gamma$ is torsion free then the setup is a precise
combinatorial analogue of that in Section 4.  For then
$\Gamma\bs\Delta$ is a finite cell complex of nonpositive curvature
with universal covering space $\Delta$ and fundamental group $\Gamma$.
Moreover $\Gamma_0\bs F_0$ is homeomorphic to the $r$-torus, and {\bf(B3)}
implies that the natural mapping $\Gamma_0\bs F_0\to\Gamma\bs\Delta$ is an
embedding of the $r$-torus into $\Gamma\bs\Delta$.\\
\item"{(b)}" The sole reason for assuming that $\Gamma_0$ is abelian in
condition {\bf(B2)} is to obtain an abelian von Neumann algebra
$\VN(\Gamma_0)$.
Everything else works equally well without this assumption. \\
\item"{(c)}" Let $\Gamma$ be a group of automorphisms of $\Delta$ which acts
properly
discontinuously and cocompactly on $\Delta$. (This is the case if
condition {\bf(B2)} is satisfied.) An apartment $F$ in $\Delta$ is called
{\it periodic} (or {\it $\Gamma$-closed}) if
the stabilizer $\Gamma_F$ of $F$ acts cocompactly on $F$.  The existence
of an abundance of periodic apartments follows from [1, Theorem 8.9]. If
$F$ is a
periodic apartment then the group $\Gamma_F$, being a Bieberbach group,
necessarily
contains a finite index subgroup $\Gamma_0\cong \ZZ^r$, which also acts
cocompactly
on $F$. Condition {\bf (B2)} is therefore satisfied for many apartments $F_0$
and
subgroups $\Gamma_0<\Gamma$.
\endroster

\endremark

The assumption {\bf(B3)} has a simple interpretation in terms of the action of
$\Gamma$ on $\Delta$.

\proclaim{Lemma 5.2}
Condition {\bf(B3)} is equivalent to the following statement.
\roster
\item"{\bf(I)}" If $\gamma\in\Gamma$, $a\in F_0^0$ and $\gamma a\in F_0^0$,
then $\gamma\in\Gamma_0$.
\endroster
\endproclaim

\demo{Proof}
Assuming {\bf(B3)}, let $\gamma\in\Gamma$, $a\in F_0^0$ and $\gamma a\in
F_0^0$.
Then $\Gamma a=\Gamma\gamma a$, and so $\Gamma_0 a=\Gamma_0\gamma a$, by
injectivity.
Thus $\gamma a =\gamma_0 a$ for some $\gamma_0\in\Gamma_0$.  However, $\Gamma$
acts
freely on $\Delta^0$.  Therefore $\gamma=\gamma_0\in\Gamma_0$.

Conversely, if {\bf(I)} holds, suppose that $a$, $a'\in F_0^0$ and $\Gamma
a=\Gamma a'$.
Then $\gamma a=a '$, for some $\gamma\in\Gamma$.  In particular, $\gamma a\in
F^0_0$.
Therefore $\gamma\in\Gamma_0$, and so $\Gamma_0 a=\Gamma_0 a'$.
\enddemo

Our next aim is to  give a combinatorial analogue of Theorem 4.9. We
begin with
some preliminary results.

\proclaim{Lemma 5.3  ([17, Lemma 2.2])}
Let $C>0$ and let $S$, $S'$ be sectors in $\Delta$.  Then either $S$ and $S'$
contain a common subsector or $S$ contains a subsector all of whose points are
at a distance greater than $C$ from $S'$.
\endproclaim

\demo{Proof}
Choose subsectors $S_1$ and~$S_1'$ of $S$ and~$S'$
respectively which lie in a common apartment [18, Chapter 9, Proposition
(9.5)].
If $S_1$ and~$S_1'$ point in the same direction, then they have a
common subsector, which is also a common subsector of~$S$ and~$S'$.

Otherwise, fix a finite $C_1>0$ so that $d(v,S_1')\leq C_1$ for
any~$v\in S'$ [18, Chapter 9, Lemma (9.2)].  Choose a subsector $S_2$
of~$S_1$ all of whose points are at a distance greater than $C+C_1$ from
$S_1'$.
Then those points are all at a distance greater than $C$
from $S'$.
\enddemo

\proclaim{Lemma 5.4 (c.f. Lemma 4.7)}
Let $F$, $F_1,\ldots,F_m$ be apartments in $\Delta$ such that, for some
$\delta>0$,
$F^0\displaystyle\operatornamewithlimits{\subset}_{\delta} F_1^0\cup\cdots\cup
F_m^0$.  If
$S$ is a sector
in
$F$ then there exists a subsector $S^*\subset S$ such that $S^*\subset F_j$,
for some $j$.
\endproclaim

\demo{Proof}
For $1\leq j\leq m$, express $F_j$ as a finite union of sectors.  Let
$\{S_\alpha:\alpha\in I\}$
denote the set of all such sectors.

Suppose that the sector $S$ does {\it not} contain a subsector in common with
any $S_\alpha$.  By Lemma 5.3, for each $\alpha\in I$ there exists
a subsector $S^*_\alpha$ of $S$, all of whose points are at distance $>\delta$
from $S_\alpha$.

Now $T=\cap_{\alpha\in I}S_\alpha^*$ is a (nonempty) subsector of $S$.
Choose a vertex $t\in T^0$.  Then $d(t,F_j)>\delta$ for each $j$.
This contradicts the assumption that
$F^0\displaystyle\operatornamewithlimits{\subset}_{\delta}F_1^0\cup\cdots\cup
F_m^0$.
\enddemo

\proclaim{Corollary 5.5}
If $F_1$, $F_2$ are apartments in $\Delta$ and
$F_1\displaystyle\operatornamewithlimits{\subset}_{\delta}
F_2$,
for some $\delta>0$, then $F_1=F_2$.
\endproclaim

\demo{Proof}
Express $F_1$ as a finite union of sectors $\{S_\alpha:\alpha\in I\}$, based
at a vertex $v\in\Delta$.  By Lemma 5.4, each $S_\alpha$ contains a
subsector $S_\alpha^*\subset F_2$.  In particular $F_1\cap F_2\neq\emptyset$
and we may assume from the start that $v\in F_1\cap F_2$.

Now for each $\alpha\in I$, $F_2$ contains $v$ and $S_\alpha^*$, and hence
also $S_\alpha$, which is the convex hull of $v$ and $S_\alpha^*$.  Thus
$F_2\supseteq F_1$.
However $F_1$, $F_2$ are isomorphic Coxeter complexes in $\Delta$.  Therefore
$F_1=F_2$.
\enddemo

Before proceeding, recall that $VN(\Gamma)$ is a $\tto$ factor if and only
if the group $\Gamma$ is ${\hbox{I.C.C.}}$.  If $\Gamma$ were a lattice in a
$p$--adic
Lie
group then the argument of [9, Lemma 3.3.1] (which uses the Borel
density theorem) could be modified to prove that $\Gamma$ is
${\hbox{I.C.C.}}$.
However not
all the groups considered in this section are embedded in a natural way as
subgroups of $p$-adic linear groups.  We therefore use a geometric argument
to verify the ${\hbox{I.C.C.}}$ property of $\Gamma$.

\proclaim{Lemma 5.6}
Let $\Delta$ be a euclidean building.  Let $\Gamma$ be a group of
automorphisms
of $\Delta$ which acts cocompactly on $\Delta$.  Then $\Gamma$ is I.C.C.
\endproclaim

\demo{Proof}
We have $\Gamma \cK = \Delta$, where $\cK\subset \Delta$ is compact.
Let $x\in \Gamma - \{e\}$, and suppose that $C= \{y^{-1}xy : y\in \Gamma\}$ is
finite.
Let
$$
\delta=\max\{d(\kappa, y^{-1}xy\kappa) : \kappa\in \cK, y\in \Gamma\}.\tag{5.1}
$$
Then
$$
d(y\kappa, xy\kappa) = d(\kappa, y^{-1}xy\kappa) \le \delta,\tag{5.2}
$$
 for all $y\in \Gamma, \kappa\in \cK$.
Therefore, for all $\xi\in \Delta$,
$$
d(\xi, x\xi) \le \delta. \tag{5.3}
$$

Choose $\eta\in \Delta$ such that $x\eta\ne \eta$
and choose an apartment $F$ in $\Delta$ with $\eta\in F$, $x\eta\notin F$.
Now by  (5.3), $F\displaystyle\operatornamewithlimits{\subset}_{\delta}xF$.
Corollary 5.5
therefore
implies that $F= xF$. In particular $x\eta\in F$, a contradiction.
\enddemo

\remark{Remark \rom{5.7}}
The proof of Lemma 5.6 also applies to a cocompact group $\Gamma$ of
isometries of a
symmetric space. (The analogue of Corollary 5.5 is [13, Lemma
5.4].) In particular
one obtains a proof of [9, Lemma 3.3.1] in the cocompact case which
avoids the use of
Borel's density theorem.
\endremark

\proclaim{Theorem 5.8}
Let $\Gamma$ be a group of automorphisms of a locally finite euclidean
building $\Delta$.
Assume that {\bf (B1), (B2), (B3)} hold.  Then $VN(\Gamma_0)$ is a strongly
singular masa of the
$\tto$ factor $VN(\Gamma)$.
\endproclaim

\demo{Proof}
In view of Lemma 4.6,
it suffices to verify condition {\bf(C2)}.
Suppose that $x_1,\ldots,x_m\in\Gamma$ and $\delta>0$ satisfy
$F_0^0\displaystyle\operatornamewithlimits{\subset}_{\delta} x_1F_0^0\cup
x_2F_0^0\cup\cdots\cup
x_mF_0^0$.
Let $S$ be a sector in $F_0$.  By Lemma 5.4, there exists a subsector
in
$S^*\subset S$ such that $S^*\subset x_jF_0$, for some $j$.  Now $S^*\subset
F_0\cap x_jF_0$.
Choose a vertex $s$ of $S^*$.  Then $s=x_j a$ for some vertex $a\in F_0^0$.
In particular $a\in F_0^0$ and $x_j a\in F_0^0$.
It follows from Lemma 5.2 that $x_j\in \Gamma_0$.
Therefore condition {\bf(C2)} is satisfied.
\enddemo

\example{Example 5.9}
{\it Groups acting on buildings of type $\tilde
A_2$}.

Suppose that the building $\Delta$ has the property that there is a group
$\Gamma$ of
automorphisms of $\Delta$ which acts freely and transitively on the vertex set
$\Delta^0$.
For buildings of type $\tilde A_2$, groups with this property have been
intensively
studied in [5, I,II].  Suppose in addition that $\Gamma$ has an
abelian
subgroup $\Gamma_0$ which acts transitively on the vertex set of an apartment
in $\Delta$.
Then {\bf(B1), (B2), (B3)} hold and so, by Theorem 5.8,
$\VN(\Gamma_0)$
is a strongly singular masa of $\VN(\Gamma)$.  Of the groups acting on
$\tilde A_2$ buildings which are enumerated in [5, II], those
labeled (4.1), (5.1), (6.1), (9.2), (13.1), (28.1) in that article
contain such a subgroup $\Gamma_0$. The groups (4.1), (5.1), (6.1)
are lattices in $\pgl_3(\QQ_3)$, but the groups (9.2), (13.1), (28.1)
do not have a natural embedding into a linear group. These groups all
have 3-torsion, but act freely on the vertex set of $\Delta$.
\endexample


\example{Example 5.10}
{\it A strongly singular but not ultrasingular masa in a $\tto$
factor with property (T)}.

Let $\Gamma$ be the group denoted by (4.1) in [5, II].
Then $\Gamma$ is a lattice subgroup of $\pgl_3(\QQ_3)$, and acts
freely and transitively on the vertex set of a building $\Delta$
of type $\tilde A_2$.  The presentation of $\Gamma$ given in [5, II]
has 13 generators $x_i$, $(0\leq i\leq 12)$, and 20 relations, among which are
$$\matrix
x_2x_4x_7&=&1&=&x_2x_7x_4\,,\\
x_2x_3x_9&=&1&=&x_2x_9x_3\,,\\
x_7x_{11}x_{12}&=&1&=&x_7x_{12}x_{11}\,,\\
x_6x_8x_9&=&1&=&x_6x_9x_8\,.
\endmatrix$$
The abelian subgroups $\Gamma_1=\langle x_2,x_4,x_7\rangle$,
$\Gamma_2=\langle x_2,x_3,x_9\rangle$, $\Gamma_3=\langle
x_7,x_{11}x_{12}\rangle$,
$\Gamma_0=\langle x_6,x_8,x_9\rangle$ are all free abelian of rank 2 and each
acts
transitively on the vertex set of an apartment.  Thus each $\VN(\Gamma_i)$ is
a
strongly singular masa of $\VN(\Gamma)$.  (Of course this implies that each
$\Gamma_i$
is a maximal abelian subgroup of $\Gamma$.)

According to [5, II, \S 5], the group $\Gamma$ has an automorphism
$f:x_i\mapsto x_{\pi(i)}$, where $\pi$ is the permutation
$(0\,\,1\,\,10\,\,5)(6\,\,11\,\,8\,\,12)(3\,\,4)(7\,\,9)$.  The action of $f$
on $\Gamma$ interchanges $\Gamma_1$ and $\Gamma_2$, and
$f^2|_{\Gamma_1}=\hbox{{\it id}}$, $f^2|_{\Gamma_2}=\hbox{{\it id}}$.
Moreover $f$ interchanges $\Gamma_3$ and $\Gamma_0$ and $f^2|_{\Gamma_0}$
exchanges the generators $x_6$ and $x_8$.

We claim that $f^2$ is an outer automorphism of $\Gamma$.  For suppose that
$f^2x=gxg^{-1}$ where $g\in\Gamma$.  Since $f^2|_{\Gamma_1}=\hbox{{\it id}}$,
and $\Gamma_1$ is a maximal abelian subgroup of $\Gamma$, we must have
$g\in\Gamma_1$.  Similarly $g\in\Gamma_2$, since $f^2|_{\Gamma_2}=\hbox{{\it
id}}$.  Thus $g\in\Gamma_1\cap\Gamma_2=\{1\}$.

It follows that $f^2$ induces an outer automorphism $\alpha$ of $\VN(\Gamma)$,
under which $\VN(\Gamma_0)$ is invariant.  In particular $\VN(\Gamma_0)$ is a
strongly singular masa of $\VN(\Gamma)$ which is {\it not} ultrasingular in
the sense of [16].  This answers in the negative a question raised in
[17, Remark 2.10].

Note that as a lattice in a higher rank group, $\Gamma$ has Kazhdan's property
(T), so that $\VN(\Gamma)$ does contain ultrasingular masas by [16, Corollary
4.5].
\endexample

In the examples above the group $\Gamma$ has torsion. It is worth examining
some cases where $\Gamma$ is torsion free.

\example{Example 5.11}
{\it A torsion free lattice in $\pgl_3$}.

Let $\G$ be the Regular $\tilde A_2$ group, which is a
lattice subgroup of $\PGL(3,\KK)$, where $\KK$ is the Laurent series field
${\bold F}_4((X))$ over the field ${\bold F}_4$ with four elements. This
group is
described in [5, I, Section 4]
and the embedding of $\G$ in $\PGL(3,{\bold F}_4((X)))$ is essentially
unique,
by the Strong Rigidity Theorem of Margulis.
The group $\G$ is torsion free and has 21 generators
$x_i, 0\le i \le 20$, and relations (written modulo 21):
$$
\cases
x_jx_{j+7}x_{j+14}=x_jx_{j+14}x_{j+7}=1 \quad & 0\le j\le 6,\\
x_jx_{j+3}x_{j-6}=1 \quad & 0\le j\le 20.\endcases
$$
This group $\G$ acts freely and transitively on the vertex set of its building
$\D$.

It follows from the first seven pairs of relations above that,
for each $j$ with $0\le j\le 6$, the generators $x_j, x_{j+7}, x_{j+14}$
pairwise commute and generate a free abelian subgroup of rank two inside $\G$,
satisfying the hypotheses of Theorem 5.8.
\endexample

\example{Example 5.12} {\it Groups acting on products of trees}.

Consider some specific examples studied in [14].  In [14, Section
3],
there is constructed a lattice subgroup $\G$ of $G = PGL_2(\QQ_p) \times
PGL_2(\QQ_l)$, where $p,l \equiv 1 \pmod 4$ are two distinct
primes. This restriction is made because $-1$ has a square root in $\QQ_p$
if and only if $p \equiv 1 \pmod 4$. The building $\D$ of $G$ is a product
of two homogeneous trees
$T_1$, $T_2$ of degrees $(p+1)$ and $(l+1)$ respectively (that is, a euclidean

building $\D$ of type $\tilde A_1\times \tilde A_1$) and $G$ is a subgroup of
$\aut(\D)$. The group $\G$ is a torsion free group which acts freely and
transitively on the vertex set $\D^0$, but which is {\it not} a product of
free groups. In fact it is an irreducible lattice in $G$.

Here is how $\G$ is constructed [14]. Let
$\HH(\ZZ)=\{\alpha=a_0+a_1i+a_2j+a_3k ; a_j\in \ZZ\}$,
the ring of integer quaternions. Let $i_p$ be a square root of $-1$ in $\QQ_p$
and define
$$\psi : \HH(\ZZ) \to PGL_2(\QQ_p) \times PGL_2(\QQ_l)$$
 by
$$\psi(a_0+a_1i+a_2j+a_3k)=\left(
\matrix
a_0+a_1i_p & a_2+a_3i_p \\
-a_2+a_3i_p  & a_0-a_1i_p
\endmatrix,
\matrix
a_0+a_1i_l & a_2+a_3i_l \\
-a_2+a_3i_l  & a_0-a_1i_l
\endmatrix
\right)
$$

Let $\tilde\G=\{\alpha=a_0+a_1i+a_2j+a_3k\in \HH(\ZZ) ; a_0\equiv 1 \pmod 2,
a_j\equiv 0 \pmod 2, j=1,2,3, |\alpha|^2=p^rl^s\}$.
Then $\G=\psi(\tilde\G$) is a torsion free cocompact lattice in $G$.
Let

$A=\{a=a_0+a_1i+a_2j+a_3k\in \tilde\G ; a_0>0, |a|^2=p\}$,

$B=\{b=b_0+b_1i+b_2j+b_3k\in \tilde\G ; b_0>0, |b|^2=l\}$.

Then $A$ contains $p+1$ elements and $B$ contains $l+1$ elements.
The images $\underline A$, $\underline B$ of $A, B$ in $\G$ generate
free groups $\G_p$, $\G_l$ of orders $p+1$, $l+1$ respectively and $\G$
itself is generated by $\underline A \cup \underline B$.
The 1-skeleton of $\D$ is the Cayley graph of $\G$ relative to this set
of generators.

By abuse of notation, identify a quaternion in $\tilde \G$ with its image in
$\G$.
(If one quaternion is a rational multiple of the other then they have the
same image in $\G$.)

It is now easy to exhibit copies of $\G_0\cong \ZZ^2$ inside $\G$, with $\G_0$
acting
freely and transitively on the vertex set of an apartment in $\D$, and
therefore
satisfying the hypotheses of the Theorem 5.8.  There are integers
$a_0, a_1$
(essentially unique) with $a_0$ odd, $a_1$ even and $a_0^2+a_1^2=p$. (The Two
Square Theorem.)
Similarly there are $b_0, b_1$ with $b_0^2+b_1^2=l$.
Let $a=a_0+a_1 i$, $b=b_0+b_1 i$. Then we can take $\G_0=\langle a, b\rangle$.
There is nothing special about the choice of $i$ rather than $j$ or $k$. Thus
we get two other possible groups $\G_0$.
Specific Example: $p=5, l=13$, $a=1+2i$, $b=3+2i$.
\endexample

\head 6. Borel subgroups of linear algebraic groups\endhead

J. ~Dixmier [6] constructed examples of singular masas by considering
groups of
homographies.  The purpose of this section is to extend his construction.  A
basic
example is the following.

\example{Example 6.1}
Let $\Gamma$ be the upper triangular subgroup of $\psl_n(\QQ)$, $n\geq 2$, and
let $\Gamma_0$ be the diagonal subgroup of $\Gamma$.  Then $\VN(\Gamma_0)$ is
a strongly singular masa of the $\tto$ factor $\VN(\Gamma)$.
[Diximer deals with the case $n=2$.]

We shall prove this result by using Proposition 4.2,
and the methods of the previous section.  In order to do this,
we let $\Gamma$ act on an appropriate euclidean building.  Choose
a prime $p$ and let $G=\psl_n(\QQ_p)$.  Then $G$ acts upon its
Bruhat-Tits building $\Delta$, whose vertex set is $G/K$, where
$K=\psl_n(\ZZ_p)$.
Here $\QQ_p$ is the $p$-adic field and $\ZZ_p$ the $p$-adic integers. Details
can
be found in [3, VI.9F].

Choose the apartment $F_0$ of $\Delta$ whose vertices are all the cosets of
the form
$$\left[\matrix p^{j_1}&&&\\
                  &p^{j_2}&&\\
                &&\ddots&\\
                &&&p^{j_n}\endmatrix \right]\cdot K,
        \qquad j_k\in\ZZ,\,1\leq k\leq n\,.$$
Then $\Gamma_0$ clearly acts transitively on the vertex set $F_0^0$.

The boundary $\Omega$ of $\Delta$ is the quotient space $G/B$, where $B$ is
the Borel subgroup of upper triangular matrices in $G$.  It is important to
note that since $\Gamma$ is a subgroup of $B$, there is a boundary point
$\omega_0$
(the coset $B$) which is stabilized by $\Gamma$.
\endexample

Consider now the following general setup.  Let $\Delta$ be a euclidean
building and let $G$ be a strongly transitive type preserving subgroup of
$\aut(\Delta)$ [7, \S 17]. This means that $G$ acts transitively on
the set of pairs $(C,F)$ where $F$ is an apartment of $\Delta$ and $C$ is a
chamber contained in $\Delta$.  Fix an apartment $F_0\subset\Delta$ and a
sector
$S_0\subset F_0$.  Then $S_0$ represents a boundary point $\omega_0$ of
$\Delta$.

Consider the {\it Borel subgroup} $B=\{g\in G:g\omega_0=\omega_0\}$.  Let
$N=\{g\in G:g F_0=F_0\}$ and let ${\frak A}=B\cap N$, the {\it Cartan
subgroup}.
Then by [7, Theorem 17.3],
$${\frak A}=\{g\in G; \hbox{$g$ acts on $F_0$ by translations}\}\,.$$
Example 6.1 is a special case of this setup, with $G=\psl_n(\QQ_p)$
[7, \S 19].

\proclaim{Theorem 6.2}
Under the above assumptions, let $\Gamma$ be an I.C.C. subgroup of $B$ and let
$\Gamma_0=\Gamma\cap N$.  Suppose also that $\Gamma_0\bs F_0^0$ is finite.
Then $\VN(\Gamma_0)$ is a strongly singular masa of $\VN(\Gamma)$.
\endproclaim

\remark{Remark \rom{6.3}}
 The group $B$ itself is ${\hbox{I.C.C.}}$. One way to see this is
to note
that
$B$ acts on $\Delta$ with finitely many vertex orbits [7, Theorem
17.6],
and so $B$ acts cocompactly on $\Delta$.
The {\rm I.C.C.} property follows from Lemma 5.6.
\endremark

\remark{Remark \rom{6.4}}
 Theorem 6.2 applies in particular to Example
6.1,
with $\Gamma$ the group of upper triangular matrices in $\psl_n(\QQ)$.
In Example 6.1 $\Gamma_0\bs F_0^0$ is a singleton, since $\Gamma_0$
acts
transitively on $F_0^0$.
The action of $B$ on $\Delta$ is continuous and the vertex set of $\Delta$ is
discrete.
Therefore the dense subgroup $\Gamma$ of $B$ has finitely many vertex orbits,
since $B$
does. Thus $\Gamma$ is {\rm I.C.C.} by Lemma 5.6, since it
acts
cocompactly on $\Delta$.
\endremark

\demo{Proof}(Theorem 6.2.) We verify conditions {\bf (C1)}, {\bf
(C2)}
 of Proposition 4.2,
with $X=\Delta^0$ and $A=F_0^0$.  The fact that $\Gamma_0\bs F_0^0$ is
finite implies condition {\bf (C1)}.

Turning to {\bf (C2)}, let $x_1,\ldots,x_m\in\Gamma$, $\delta>0$ and
$F_0^0\displaystyle\operatornamewithlimits{\subset}_{\delta}x_1 F_0^0\cup
x_2F_0^0\cup\cdots\cup
x_mF_0^0$.
Choose a sector $S$ in $F_0$ opposite $S_0$.  By Lemma 5.4, there
exists a subsector $S^*\subset S$ such that $S^*\subset x_j F_0$ for some $j$.

Since $x_j\omega_0=\omega_0$, the two sectors  $x_j S_0$ and $S_0$ have a
common subsector $S^*_0$.

We now have $S^*_0\cup S^*\subset F_0$ and $S_0^*\cup S^*\subset x_j F_0$.
However, opposite sectors in an apartment determine that apartment
completely [3, VI.9, Lemma 2 and IV.5, Theorem 1].  Therefore
$x_jF_0=F_0$.
In other words, $x_j\in\Gamma\cap N=\Gamma_0$.  This establishes condition
{\bf (C2)}.
\enddemo

\noindent {\it Remark.}
(a) Generalizing Example 6.1, one can clearly let $\Gamma$ be the
upper triangular subgroup of $\psl_n(\KK)$, where $\KK$ is any subfield of
$\QQ_p$
for some prime $p$, with $\QQ\subset\KK\subset\QQ_p$.  In fact $\KK$ could be
an
appropriate subfield of any nonarchimedean local field.  The abelian subgroup
$\Gamma_0$
is again the diagonal subgroup of $\Gamma$.
Note that $\Gamma$ is amenable and so $\VN(\Gamma)$
is the hyperfinite $\tto$ factor.\newline
\noindent (b) Other generalizations are possible.  For example, one could
replace $\QQ_p$
by $\RR$, and work with symmetric spaces.  In the case $n=2$, $\Gamma$ would
be the
subgroup of upper triangular matrices in $\psl_2(\QQ)$, and $\Gamma_0$ its
diagonal subgroup.
These groups act on the hyperbolic plane $X$ and the crucial point in the
argument
is the fact that a geodesic in $X$ is uniquely determined by its two boundary
points.

\enddocument